\documentclass[12pt]{article}
\usepackage{e-jc+}

\title{The Largest Possible Finite Degree of Functions between Commutative Groups}

\author{Uwe Schauz\\
\small Department of Pure Mathematics\\[-0.8ex]
\small Xi’an Jiaotong-Liverpool University\\[-0.8ex]
\small Suzhou 215123, China\\[-0.8ex]
\small \texttt{uwe.schauz@xjtlu.edu.cn}}

\date{\printday}

\usepackage{fixltx2e}

\usepackage[ansinew]{inputenc}
\usepackage{srcltx}
\usepackage[T1]{fontenc}
\usepackage{lmodern}
\usepackage{exscale}

\usepackage{ngerman}
\usepackage{textcomp}

\usepackage{amsmath}
\usepackage{amsthm}
\usepackage{extpfeil}

\usepackage{amssymb}
\usepackage{stmaryrd}
\usepackage{bm}

\usepackage{enumerate}

\usepackage{afterpage}
\usepackage{xspace}

\newtheoremstyle{Theorem}{.7\baselineskip}{1\baselineskip}{\itshape}{}{\bfseries}{.}{ }{}
\theoremstyle{Theorem}
\newtheorem{Satz}{Theorem}[section]

\newtheorem{Korollar}[Satz]{Corollary}

\newtheorem{Lemma}[Satz]{Lemma}

\newtheoremstyle{Definition}{.6\baselineskip}{.8\baselineskip}{}{}{\bfseries}{.}{ }{}
\theoremstyle{Definition}
\newtheorem{Definition}[Satz]{Definition}
\theoremstyle{definition}
\newtheorem{Beispiel}[Satz]{Example}

\theoremstyle{remark}

\newenvironment{sequation*}{\begin{small}\begin{equation*}}{\end{equation*}\end{small}}
\newenvironment{Beweis}[1][Proof]{\begin{proof}[#1]}{\end{proof}}

\newenvironment{proofsize}[1]{\begin{small}#1}{\end{small}}

\newcommand\taghere{\stepcounter{equation}\tag{\theequation}}

\newcommand\´{\kern 1pt}
\newcommand\ssp{\kern 1pt}

\newcommand\noms{\hspace{-\mathsurround}}

\newcommand\printday{\today\xspace}
\newcommand\Rand[1]{
  \marginpar{\raggedleft\scriptsize\hspace{0pt}#1}}%

\renewcommand{\(}{\noms$}
\renewcommand{\)}{\noms$}
\renewcommand\frac[2]{\genfrac{}{}{.4pt}{}{#1}{#2}}

\renewcommand\dfrac[2]{\genfrac{}{}{.4pt}{0}{#1}{#2}}

\newcommand\mathRand[1]{\hspace{\mathsurround}\Rand{#1}\nolinebreak\noms}
\def\rand #1"#2"{\mathRand{\(#2\)}#1#2}
\def\randd #1"#2"#3\randd#4"#5"{\mathRand{\(#2\), \(#5\)}#1#2#3#4#5}

\newcommand\eqby[2][=]%
  {\ensuremath{\overset{\makebox[0pt]{\ensuremath{\smash[t]{\scriptstyle#2}}}}{#1}}}

\def\eqbelow#1#2{\underset{\overset{\shortparallel}{#2}}{\smash{#1}}}

\newcommand\Z{\mathbb{Z}}
\newcommand\N{\mathbb{N}}

\renewcommand\sb{\subseteq}

\newcommand\mto{\mapsto}
\newcommand\lmto{\longmapsto}
\newcommand\lto{\longrightarrow}

\newcommand\lEqi{\Longleftrightarrow}
\newcommand\nach{\mathbin\circ}

\newcommand\fa{\forall\,}
\newcommand\ex{\exists\,}

\newcommand\DP{\colon\discretionary{\!\kern -.17em}{}{}}
\newcommand\mitsymbol{\textup{\textbrokenbar}}
\renewcommand\mit{\,\ \discretionary{\mitsymbol}{}{}\mitsymbol\ \,}

\renewcommand\div{\mathrel{\bigm\lfloor\!\!\!\bigm\lfloor}}
\newcommand\vid{\mathrel{\bigm\rfloor\!\!\!\bigm\rfloor}}
\newcommand\ndiv{\mathrel{\;\!\div\hspace{-12pt}\kern0pt\lower2pt%
  \hbox{\ensuremath{^\diagup}}\!}}
\newcommand\ndivps{\mathrel{\;\!\div\hspace{-9pt}\kern0pt\lower2pt%
  \hbox{\ensuremath{^\diagup}}\!}}
\newcommand\nvid{\mathrel{\;\!\vid\hspace{-12pt}\kern0pt\lower2pt%
  \hbox{\ensuremath{^\diagup}}\!}}
\newcommand\nvidps{\mathrel{\;\!\vid\hspace{-9pt}\kern0pt\lower2pt%
  \hbox{\ensuremath{^\diagup}}\!}}

\providecommand\abs[1]{\lvert#1\rvert}
\providecommand\Abs[1]{\bigl\lvert#1\bigr\rvert}

\DeclareMathOperator\fdeg{\textsc{Fdeg}}
\DeclareMathOperator\pdeg{\textsc{Pdeg}}

\newcommand\p{\varphi}

\renewcommand\b{\beta}
\newcommand\ä{\alpha}

\renewcommand\d{\delta}

\newcommand\e{\varepsilon}

\begin{document}
\maketitle

%
%
%
%
%
\begin{abstract}
We consider maps between commutative groups and their functional degrees.
These degrees are defined based on a simple idea -- the functional degree should
decrease if a discrete derivative is taken. We show that the maps
of finite functional degree are precisely the maps that can be written as polyfracts, as
polynomials in several variables but with binomial functions in the place of powers.
Moreover, the degree of a
polyfract coincides with its functional degree. 
We use this to determine the largest possible finite functional degree that the maps
between two given finite commutative groups can have. This also yields a solution to
Aichinger and Moosbauer's problem of finding the nilpotency degree of the
augmentation ideal of the group ring $\Z_{p^\b}[\Z_{p^{\ä_1}}\!\times\Z_{p^{\ä_2}}
\!\times\dotsm\times\Z_{p^{\ä_n}}]$. 
Some generalizations and simplifications of proofs to underlying facts are presented,
too.
\end{abstract}
%
%
%

\section{Introduction}\label{sec.int}

In the recent paper \cite{aimo}, Aichinger and Moosbauer studied what they call the
\emph{functional degree} $\fdeg(f)$ of maps $f\DP A\lto B$ between commutative
groups $(A,+)$ and $(B,+)$, in order to generalize some Chevalley-Warning type
results. The functional degree was examined in other papers before, e.g.\ in \cite{la,ma},
partially in less generality. Aichinger and Moosbauer, however, introduced the name
functional degree (which reflects the fact that its definition does not rely on a term
representation of $f$) and investigated its properties more comprehensively. In
\cite[Problem\,8.3]{aimo}, they also raised the question of how large the functions
degree possible can be, if we restrict ourselves to functions between to given finite
commutative \(p\)-groups. 
Equivalently, one can also investigate the nilpotency degree of the augmentation ideal
of the group ring $\Z_{p^\b}[\Z_{p^{\ä_1}}\!\times\Z_{p^{\ä_2}}
\!\times\dotsm\times\Z_{p^{\ä_n}}]$ ($p$ prime), as explained in \cite[Lem.\,7.3]{aimo},
but we do not look at the problem from this side. We generalize Aichinger and
Moosbauer's question and ask for the largest possible finite functional degree of
functions between arbitrary finite commutative groups. Without the restriction to
\(p\)-groups, the functional degree may become infinite, but there must be a largest
degree among all finite degrees (because there are only finitely many functions between
two given finite groups). We ask what this upper limit is, and answer this question in our
last section, Section\,\ref{sec.max}. We do this right after presenting and explaining the
classification of all functions of finite degree in Section\,\ref{sec.clas}, according to
\cite{aimo} and \cite{sch}. All that is based on a term representation of functions of finite
functional degree, which we provide in Section\,\ref{sec.trep}. There we learn how to
write functions of finite functional degree as polynomials in several variables but with
binomial functions $\tbinom{x_j}{\d_j}$ in the place of powers $\smash{x_j^{\d_j}}$\!. We
call this kind of polynomials \emph{polyfracts} and cite many results about them from
\cite{sch}. We also present some generalizations of results in \cite{sch}, sometimes with
different proofs. We start our investigations with the introduction of the functional degree
in the next section, Section\,\ref{sec.fdeg}.

\section{The Functional Degree}\label{sec.fdeg}

In this section, we introduce the functional degree $\fdeg(f)$ of functions $f\DP A\lto B$
between commutative groups $(A,+)$ and $(B,+)$. The definition of the functional
degree in \cite[Def.\,2.1]{aimo} is given on an abstract level, in the language of group
rings, modules and ideals. It can, however, also be stated with the help of
\emph{(discrete) derivatives} and \emph{difference operators}. For each $g\in A$, the
difference operator \rand$"\Delta_g"\DP B^A\lto B^A$, $f\lmto\Delta_gf$ is defined via
\begin{equation}
[\Delta_g f](x)\,:=\,f(x+g)-f(x)\,\ .
\end{equation}
With this notation, we can reinterpret \cite[Lem.\,2.2]{aimo} as a definition:

\begin{Definition}\label{def.1} 
Let $A$ and $B$ be commutative groups and $G$ a generating subset of $A$. The
\emph{functional degree} \rand$"\fdeg(f)"$ of a map $f\DP A\lto B$ is the smallest
number $m\geq0$ such that, for all $g_1,g_2,\dotsc,g_{m+1}\in G$,
$$\Delta_{g_1}\Delta_{g_2}\dotsm\Delta_{g_{m+1}} f\,\equiv\,0\,,$$
if such an $m$ exists, and $\fdeg(f)=\infty$ otherwise. In particular, we have
$\fdeg(0)=0$.
\end{Definition}

It follows from \cite[Lem.\,2.2]{aimo} that the functional degree does not depend on the
choice of the set of generators $G$. We will use that to make a special choice for $G$
in our main definition below, for finitely generated commutative groups $A$. Up to
isomorphy, such groups can be written as
\begin{equation}
A=\Z_{q_1}\!\times\Z_{q_2}\!\times\dotsm\times\Z_{q_n}\qquad\text{with integers $q_i\geq0$,}
\end{equation}
where we usually allow the extremal cases $q_i=1$ and $q_i=0$, i.e.\
  \rand\rand\rand\begin{equation}
"\Z_{q_i}":=\Z/q_i\Z\quad\text{with}\quad
"\Z_1"=\Z/\Z=\{0\}\quad\text{and}\quad"\Z_0"=\Z/\{0\}=\Z\,\ .
\end{equation}
This representation of $A$, however, is usually not uniquely determined. Choosing and
fixing one representation, as we will do it, is like choosing coordinates. It will lead us to a
coordinate dependent definition of the functional degree. With that coordinate
dependent definition, we will see that $\fdeg(f)$ actually may be seen as a kind of total
degree, with some kind of monomials and a term representation of $f$\!. To get there,
we observe that the tuples
\begin{equation}
e_i:=(0,\dotsc,0,1,0,\dotsc,0)
=\bigl([0]_{q_1},\dotsc,[0]_{q_{i-1}},[1]_{q_i},[0]_{q_{i+1}},\dotsc,[0]_{q_n}\bigr)
\end{equation}
form a generating set $G$ of
$A=\Z_{q_1}\!\times\Z_{q_2}\!\times\dotsm\times\Z_{q_n}$. Here, one actually could
omit those $e_i$ for which $q_i=1$, as they are zero, since $1=1+q_i\Z=0+q_i\Z=0$ if
$q_i=1$. One could also omit the entire factors $\Z_{q_i}$ with $q_i=1$ in the
representation $\Z_{q_1}\!\times\Z_{q_2}\!\times\dotsm\times\Z_{q_n}$ of $A$. But,
allowing such dummy coordinates will become helpful later, and the corresponding
$e_i$ do not disturb in the following definitions either (exactly as $0\in G$ is allowed and
does not disturb in Definition\,\ref{def.1}). If we abbreviate
\begin{equation}
  \underset{^i}\Delta:=\Delta_{e_i}\quad\text{i.e.}\quad[\underset{^i}\Delta f](x)\,:=\,f(x+e_i)-f(x)\,\ ,
\end{equation}
we obtain the following special case of the definition above:

\begin{Definition}
Assume $A=\Z_{q_1}\!\times\Z_{q_2}\!\times\dotsm\times\Z_{q_n}$ with integers
$q_i\geq0$. Also assume $B$ is a commutative group, and $f\DP A\lto B$ is a map. The
\emph{functional degree} $\fdeg(f)$ of $f$ is the smallest number $m\geq0$ such that,
for all $d_1,d_2,\dotsc,d_n\geq0$ with $\sum_{i=1}^nd_i>m$,
$$\underset{^1}\Delta^{\!d_1}\underset{^2}\Delta^{\!d_2}\dotsm\underset{^n}
    \Delta^{\!d_n}f\,\equiv\,0\,,$$
if such an $m$ exists, and $\fdeg(f)=\infty$ otherwise. In particular, we have
$\fdeg(0)=0$.
\end{Definition}

In \cite[Def.\,5.1]{aimo}, Aichinger and Moosbauer also define the $i^{\text{th}}$
\emph{partial functional degree} $\pdeg_i(f)$ of maps $f\DP A=\prod_{j=1}^n A_j\lto B$,
where the $A_j$ and $B$ are commutative groups.
They simply set
\begin{equation}\label{eq.pd}
\pdeg_i(f):=\sup\bigl\{\fdeg(f\nach E_a^{(i)})\!\mit\!a\in A\bigr\}
\end{equation}
where, for each fixed $a=(a_1,a_2,\dotsc,a_n)\in A$, 
$E_a^{(i)}$ is defined by
\begin{equation}
E_a^{(i)}\DP A_i\lto A_{1}\!\times A_{2}\!\times\dotsm\times A_{n}\,,\quad x\lmto(a_1,\dotsc,a_{i-1},x,a_{i+1},\dotsc,a_n)\,.
\end{equation}
If now $A_j=\Z_{q_j}$, so that
$A=\Z_{q_1}\!\times\Z_{q_2}\!\times\dotsm\times\Z_{q_n}$ as above, then for each
$m\in\N$, we have
\begin{align*}
\pdeg_i(f)\leq m&\ \lEqi\ \fa a\in A\DP\fdeg(f\nach E_a^{(i)})\leq m\\[4pt]
&\ \lEqi\ \fa a\in A\DP
  \underset{^i}\Delta^{\!m+1}(f\nach E_a^{(i)})
  \equiv0\taghere\\
&\ \lEqi\ \fa a\in A\DP(\underset{^i}\Delta^{\!m+1}f)\nach E_a^{(i)}
  \equiv0\\
&\ \lEqi\ \underset{^i}\Delta^{\!m+1}f\equiv0\,.
\end{align*}
This yields the following simple definition of the partial degree:

\begin{Definition}
Assume $A=\Z_{q_1}\!\times\Z_{q_2}\!\times\dotsm\times\Z_{q_n}$ with integers
$q_i\geq0$, $B$ is a commutative group, and $f\DP A\lto B$ is a map. The
$i^{\text{th}}$ \emph{partial functional degree} $\pdeg_i(f)$ of $f$ is the smallest number
$m$ such that,
$$\underset{^i}\Delta^{\!m+1}f\,\equiv\,0\,,$$
if such an $m$ exists, and $\pdeg_i(f)=\infty$ otherwise. In particular, we have
$\pdeg_i(0)=0$.
\end{Definition}

From this definition one can very easily deduce the inequality
\begin{equation}\label{eq.pfp}
\pdeg_i(f)\,\leq\,\fdeg(f)\,\leq\,\sum_{j=0}^n\pdeg_j(f)\,,
\end{equation}
which took two pages to prove in \cite[Th.\,5.2]{aimo}, but in the more general setting
that we mentioned in \eqref{eq.pd} above.

\section{A Term Representation}\label{sec.trep}

In this section, we introduce a term representation of functions
\begin{equation}
A=\Z_{q_1}\!\times\Z_{q_2}\!\times\dotsm\times\Z_{q_n}\lto B
\end{equation}
with finite functional degree, and show that their total and partial degrees
coincide with their functional total and partial degrees. We start with the case
$q_1=q_2=\dotsb=q_n=0$:

\subsection{Maps of the Form $A=\Z^n\lto B$}

For maps of the form $A=\Z^n\lto B$, one might first consider polynomials 
with coefficients in $B$. For this to make sense, it is enough that every
commutative group $B$ is a \(\Z\)-module, which we might regard as right
\(\Z\)-module. The substitution of integers into the symbolic variables $X_i$ of a
polynomial $\sum_{\d\in\N^n}b_\d X_1^{\d_1}X_2^{\d_2}\dotsm X_n^{\d_n}$ then
yields well-defined function values in $B$. The only problem is that we do not get
enough functions in this way. There are functions of finite functional degree that
cannot be written in this form. To improve the situation, we look at binomial
polynomials or \emph{monofracts}, as we call them. The monofract in one
variable $X$ of degree $\d\in\N$ is defined as
 \rand\begin{equation}
"\dbinom{X}{\d}"\,:=\,\dfrac{X(X-1)\dotsm(X-\d+1)}{\d!}\quad\text{where}\quad\dbinom{X}{0}\,:=\,1\,\
.
\end{equation}
For $n$ symbolic variables $X_1,X_2,\dotsc,X_n$ and for
$\d_1,\d_2,\dotsc,\d_n\in\N$\!, we set
\begin{equation}
\dbinom{X_1,X_2,\dotsc,X_n}{\d_1,\d_2,\dotsc,\d_n}
\,:=\,\dbinom{X_1}{\d_1}\dbinom{X_2}{\d_2}\dotsm\dbinom{X_n}{\d_n}\,\ .
\end{equation}
Each monofract (binomial polynomials) still gives rise to an integer valued
functions, the function
\begin{equation}
\Z^n\lto\Z\,,\quad(x_1,x_2,\dotsc,x_n)\lmto\dbinom{x_1,x_2,\dotsc,x_n}{\d_1,\d_2,\dotsc,\d_n}\,.
\end{equation}
So, it makes sense to consider linear combinations
\begin{equation}
P\,=\,\sum_{\d\in\N^n}b_\d\dbinom{X_1,X_2,\dotsc,X_n}{\d_1,\d_2,\dotsc,\d_n}\,,
\end{equation}
where only finitely many coefficients $b_\d\in B$ are non-zero. The substitution of
integers $x_i$ into the symbolic variables $X_i$ of such \emph{\(B\)-polyfracts}
$P$\!, as we call them, yields well-defined function values $P(x_1,x_2,\dotsc,x_n)$
in $B$. Hence, through substitution, we obtain to each polyfract $P$ a function
 \rand\begin{equation}
"P|_{\Z^n}"\DP\Z^n\lto B,\quad x\lmto P(x)\,.
\end{equation}
This entails that we also have a map
\begin{equation}
B\tbinom{X_1,X_2,\dotsc,X_n}{\Z^n}\lto B^{\Z^n},\quad P\lmto P|_{\Z^n}
\end{equation}
from the set of all \(B\)-polyfracts
\Rand{$B\tbinom{X_{\!1},..,X_{\!n}\!}{\Z^n}$}$B\tbinom{X_1,X_2,\dotsc,X_n}{\Z^n}$ into
the set of functions $B^{\Z^n}$\!\!. This map is even a group homomorphism between
$B\tbinom{X_1,X_2,\dotsc,X_n}{\Z^n}$ as group with coefficientwise addition and
$B^{\Z^n}$\!\! as group with pointwise addition, i.e.\ for polyfracts $P$ and $Q$,
\begin{equation}
(P+Q)|_{\Z^n}\,=\,P|_{\Z^n}+\,Q|_{\Z^n}\,.
\end{equation}
But, the map $P\lmto P|_{\Z^n}$ has another interesting property, as was shown
in \cite[Th.\,2.4]{sch} already, in less generality:

\begin{Satz}
Let $B$ be a commutative group. The homomorphism
$$
B\tbinom{X_1,X_2,\dotsc X_n}{\Z^n}\lto B^{\Z^n}\ ,\quad P\lmto P|_{\Z^n}
$$
is \emph{injective}.
\end{Satz}

\begin{Beweis}
Assuming $P\neq0$, it is enough to show that $P|_{\Z^n}\not\equiv0$. As
$P=\sum_{\d\in\N^n}b_\d\tbinom{X_1,X_2,\dotsc,X_n}{\d_1,\d_2,\dotsc,\d_n}\neq0$,
there exists an $\e\in D:=\{\d\!\mit\!b_\d\neq0\}$ that is minimal in the sense that
there is no other $\d\in D$ with $\d\leq\e$. So, every $\d\neq\e$ in $D$ is bigger
than $\e$ in at least one coordinate $i$. If we have $\d_i>\e_i$, then
\begin{equation}
\dbinom{X_1,X_2,\dotsc,X_n}{\d_1,\d_2,\dotsc,\d_n}\bigg|_{X=\e}\!=\,0\quad\text{as}\quad
\dbinom{X_i}{\d_i}\bigg|_{X_i=\e_i}
\!=\,\dbinom{\e_i}{\d_i}
\,=\,0\,.
\end{equation}
Hence, the only monofracts $\tbinom{X_1,X_2,\dotsc,X_n}{\d_1,\d_2,\dotsc,\d_n}$ in
$P$ that contributes to the function value $P(\e)$ is the one with $\d=\e$. We have
\begin{equation}
P(\e)\,=\,b_\e\dbinom{X_1,X_2,\dotsc,X_n}{\e_1,\e_2,\dotsc,\e_n}\bigg|_{X=\e}\!=\,b_\e\,\neq\,0\,,
\end{equation}
which shows that $P|_{\Z^n}\not\equiv0$, indeed.
\end{Beweis}

We see that to every map $f\in B^{\Z^n}$ there exists at most one polyfract $P\in
B\tbinom{X_1,X_2,\dotsc X_n}{\Z^n}$ that describes (interpolates) this map,
$P|_{\Z^n}=f$. So, we do not have to make a difference between polyfracts $P$ and
polyfractal maps $P|_{\Z^n}$ and will write $P$ for $P|_{\Z^n}$. We also do not have to
distinguish between the coefficientwise addition and the pointwise addition. Polyfracts
are special maps and they form a subgroup of the group of all maps from $\Z^n$ to $B$.
By identifying $P$ with $P|_{\Z^n}$, we obtain the embedding
\begin{equation}
B\tbinom{X_1,X_2,\dotsc X_n}{\Z^n}\,\sb\,B^{\Z^n}\,.
\end{equation}
In particular, we can apply the difference operators
$\smash{\underset{^i}\Delta}$ to polyfracts. By Pascal's rule we have, for
$\ell\in\N$,
\begin{equation}\label{eq.Dmi}
\underset{^i}\Delta\dbinom{X_i}{\ell+1}\,=\,\dbinom{X_i}{\ell}\quad\ \text{and}\quad\
\underset{^i}\Delta\dbinom{X_i}{0}\,=\,0\,.
\end{equation}
Using the linearity of $\smash{\underset{^i}\Delta}$, it is now easy to calculate discrete
derivatives of polyfracts. And, theoretically, this can be used to calculate the functional
degree of polyfracts. But, there is an even easier way. Polyfracts already have total and
partial degrees, exactly as normal polynomials. We have
\begin{equation}
\deg_i\tbinom{X_1,X_2,\dotsc,X_n}{\d_1,\d_2,\dotsc,\d_n}:=\d_i
\quad\text{and}\quad
\deg\tbinom{X_1,X_2,\dotsc,X_n}{\d_1,\d_2,\dotsc,\d_n}:=\d_1+d_2+\dotsb+d_n\,\,,
\end{equation}
and this is extended to polyfracts
$P=\sum_{\d\in\N^n}b_\d\tbinom{X_1,X_2,\dotsc,X_n}{\d_1,\d_2,\dotsc,\d_n}$ by taking
the maximum over all monofracts in $P$:

\begin{Definition}
For polyfracts
$P=\sum_{\d\in\N^n}b_\d\tbinom{X_1,X_2,\dotsc,X_n}{\d_1,\d_2,\dotsc,\d_n}\neq0$, we
set
\begin{equation}
\deg_i(P):=\max_{b_\d\neq0}(\d_i)
\quad\text{and}\quad
\deg(P):=\max_{b_\d\neq0}(\d_1+\d_2+\dotsb+\d_n)\,\,.
\end{equation}
We also set $\deg_i(0):=0$ and $\deg(0):=0$.\footnote{To match the conventions
regarding the zero function in \cite{aimo}, which we kept in the definitions above,
we also define $\deg(0):=0$, in contrast to the common convention of $0$ having
negative infinite degree.}
\end{Definition}

It is easy to see that this definition coincides with the definition of the functional degrees
on the set of polyfracts:
\begin{Lemma}
If $P\in B\tbinom{X_1,X_2,\dotsc X_n}{\Z^n}$, then
$$
\deg_i(P)=\pdeg_i(P)\quad\text{and}\quad\deg(P)=\fdeg(P)\,.
$$
\end{Lemma}

\begin{Beweis}
Based on Equation\,\eqref{eq.Dmi} we see that, for $d,\d\in\N^n$\!,
\begin{equation}
\underset{^1}\Delta^{\!d_1}\underset{^2}\Delta^{\!d_2}\dotsm\underset{^n}\Delta^{\!d_n}
\dbinom{X_1,X_2,\dotsc,X_n}{\d_1,\d_2,\dotsc,\d_n}\,\neq\,0
\quad\lEqi\quad d\leq\d\,.
\end{equation}
From this it follows that, for polyfracts
$P=\sum_{\d\in\N^n}b_\d\tbinom{X_1,X_2,\dotsc,X_n}{\d_1,\d_2,\dotsc,\d_n}$ and
numbers $m\in\N$, with the help of the sets $\Delta:=\{\d\in\N^n\!\mit\!b_\d\neq0\}$ and
$D:=\{d\in\N^n\!\mit\!\sum_{i=1}^nd_i>m\}$,
\begin{align*}
\fdeg(P)>m
&\ \lEqi\ \ex(d_i)\in D\DP\underset{^1}\Delta^{\!d_1}\underset{^2}\Delta^{\!d_2}\dotsm
  \underset{^n}\Delta^{\!d_n}P\neq0\\
&\ \lEqi\ \ex(d_i)\in D\DP\ex\d\in\Delta\DP d\leq\d\\[4pt]
&\ \lEqi\ \ex\d\in\Delta\DP\ex(d_i)\in D\DP d\leq\d\taghere\\[4pt]
&\ \lEqi\ \ex\d\in\Delta\DP\d\in D\\[4pt]
&\ \lEqi\ \ex\d\in\Delta\DP{\textstyle\sum_{i=1}^n\d_i>m}\\[4pt]
&\ \lEqi\ \deg(P)>m\,,
\end{align*}
so that $\fdeg(P)=\deg(P)$, and
\begin{align*}
\pdeg_i(P)>m
&\ \lEqi\ \underset{^i}\Delta^{\!m+1}P\neq0\\
&\ \lEqi\ \ex\d\in\Delta\DP m+1\leq\d_i\taghere\\[4pt]
&\ \lEqi\ \deg_i(P)>m\,,
\end{align*}
so that $\pdeg_i(P)=\deg_i(P)$.
\end{Beweis}

That, further, all functions 
of finite functional degree 
are given by polyfracts can be seen from the following kind of Taylor Theorem. It was in
less generality already stated in \cite[Th.\,2.7]{sch}, but with a mistake in the
finally printed version: 

\begin{Satz}\label{sz.tay2}
Let $d_1,d_2,\dotsc,d_n\in\N$ and $f\DP\Z^n\lto B$ be a map with $\pdeg_i(f)\leq
d_i$, for $i=1,2,\dotsc,n$, then
$$
f(x)\,\equiv\,\sum_{\d_1=0}^{d_1}\,\sum_{\d_2=0}^{d_2}\dotsb\sum_{\d_n=0}^{d_n}
[\underset{^1}\Delta^{\!\d_1}\underset{^2}\Delta^{\!\d_2}\dotsm
\underset{^n}\Delta^{\!\d_n}f](0)\,
\dbinom{x_1,x_2,\dotsc,x_n}{\d_1,\d_2,\dotsc,\d_n}\,\ .
$$
\end{Satz}

\begin{Beweis}
For $n=1$ this follows from the fact that, based on Equation\,\eqref{eq.Dmi}, for the
functions $f\DP x_1\mto f(x_1)$ and $g\DP x_1\mto
g(x_1):=\sum_{\d_1=0}^{d_1}[\Delta^{\!\d_1}f](0)\, \tbinom{x_1}{\d_1}$, we have
\begin{equation}
\Delta^{\!d_1+1}f\,\equiv\,0\,\equiv\,\Delta^{\!d_1+1}g\quad\text{and}\quad
[\Delta^{\!i}f](0)\,=\,[\Delta^{\!i}g](0)\ \,\text{for $0\leq i\leq d_1$}\,.
\end{equation}
Indeed, we just have to ``integrate'' $d_1+1$ times to deduce $f\equiv g$ from
$\Delta^{\!d_1+1}f\equiv\Delta^{\!d_1+1}g$ and the ``initial conditions''. One starts with
the ``initial condition'' $[\Delta^{\!d_1}f](0)=[\Delta^{\!d_1}g](0)$ and uses
$\Delta^{\!d_1+1}f\equiv\Delta^{\!d_1+1}g$ to show that
$[\Delta^{\!d_1}f](\pm1)=[\Delta^{\!d_1}g](\pm1)$, and then
$[\Delta^{\!d_1}f](\pm2)=[\Delta^{\!d_1}g](\pm2)$, and then
$[\Delta^{\!d_1}f](\pm3)=[\Delta^{\!d_1}g](\pm3)$, etc. Hence, one ``integration'' yields
$\Delta^{\!d_1}f\equiv\Delta^{\!d_1}g$, and $d_1$ further integrations give us $f\equiv
g$.

Now, assume we have proven the statement already for functions with up to $n-1$
variables. We can apply that induction assumtion to the function $x_n\mto
f(x_1,x_2,\dotsc,x_{n-1},x_n)$, with $x_1,x_2,\dotsc,x_{n-1}$ regarded as fixed given,
and afterwards to the functions
$(x_1,x_2,\dotsc,x_{n-1})\mto\smash{[\underset{^n}\Delta^{\!\d_n}f]_{\eqbelow{x_n}{0}}}
:=[\underset{^n}\Delta^{\!\d_n}f]\!\!\bigm|_{x_n=0}$ that arise from this step:
\begin{align*}
f(x)\equiv&\sum_{\d_n}^{d_n}[\underset{^n}\Delta^{\!\d_n}f]_{\eqbelow{x_n}{0}}\,\dbinom{x_n}{\d_n}\\
\equiv\sum_{\d_n=0}^{d_n}&\biggl[\sum_{\d_1=0}^{d_1}\!\dotsb\!\sum_{\!\!\!\d_{n-1}=0\!\!\!}^{d_{n-1}}
\Bigl[\underset{^1}\Delta^{\!\d_1}\!\dotsm
\underset{^{\!\!\!n-1\!\!\!}}\Delta^{\!\d_{n-1}}\!\bigl[[\underset{^n}
\Delta^{\!\d_n}f]_{\eqbelow{x_n}{0}}\bigr]\Bigr]_{\eqbelow{x_1}{0},\dotsc,\eqbelow{x_{n-1}}{0}}\,
\dbinom{x_1,\dotsc,x_{n-1}}{\d_1,\dotsc,\d_{n-1}}\biggr]\dbinom{x_n}{\d_n}\\
\equiv\sum_{\d_n=0}^{d_n}&\biggl[\sum_{\d_1=0}^{d_1}\!\dotsb\!\sum_{\!\!\!\d_{n-1}=0\!\!\!}^{d_{n-1}}
\Bigl[\bigl[\underset{^1}\Delta^{\!\d_1}\!\dotsm
\underset{^{\!\!\!n-1\!\!\!}}\Delta^{\!\d_{n-1}}\![\underset{^n}
\Delta^{\!\d_n}f]\bigr]_{\eqbelow{x_n}{0}}\Bigr]_{\eqbelow{x_1}{0},\dotsc,\eqbelow{x_{n-1}}{0}}\,
\dbinom{x_1,\dotsc,x_{n-1}}{\d_1,\dotsc,\d_{n-1}}\dbinom{x_n}{\d_n}\biggr]\\
\equiv\sum_{\d_1=0}^{d_1}&\,\sum_{\d_2=0}^{d_2}\dotsb\sum_{\d_n=0}^{d_n}
[\underset{^1}\Delta^{\!\d_1}\underset{^2}\Delta^{\!\d_2}\dotsm
\underset{^n}\Delta^{\!\d_n}f](0)\,
\dbinom{x_1,x_2,\dotsc,x_n}{\d_1,\d_2,\dotsc,\d_n}\,.\taghere\\[-4mm]
\end{align*}
The theorem holds for all $n$.
\end{Beweis}

As, by Inequality\,\eqref{eq.pfp}, the partial functional degrees are bounded if the functional
degree is bounded, the last theorem and the last lemma yield the following:

\begin{Satz}\label{sz.eq}
The functions $f\DP\Z^n\lto B$ with $\fdeg(f)<\infty$ are exactly the functions that are
given by \(B\)-polyfracts, i.e.\
$$
\bigl\{f\in B^{\Z^n}\!\mit\!\fdeg(f)<\infty\bigr\}\,=\,B\tbinom{X_1,X_2,\dotsc X_n}{\Z^n}\,.
$$
Moreover, if $f$ is written as polyfract, having a total degree $\deg(f)$ and partial
degrees $\deg_i(f)$, then
$$
\fdeg(f)\,=\,\deg(f)\quad\text{and}\quad\pdeg_i(f)\,=\,\deg_i(f)\quad\text{for $i=1,2,\dotsc,n$.}
$$\vspace{-4mm}
\end{Satz}

\subsection{The General Form $A=\Z_{q_1}\!\times\Z_{q_2}\!\times\dotsm\times\Z_{q_n}\lto B$}

To be able to apply the previous results, we show that the group
$B^{\Z_{q_1}\!\times\Z_{q_2}\!\times\dotsm\times\Z_{q_n}}$ with pointwise addition
may be viewed as subgroup of the group $B^{\Z^n}$ with pointwise addition. Given
the surjective group homomorphism
\begin{align*}
\p\DP\Z^n&\lto\Z_{q_1}\!\times\Z_{q_2}\!\times\dotsm\times\Z_{q_n}\taghere\\
(x_1,x_2,\dotsc,x_n)&\lmto
(x_1+q_1\Z,x_2+q_2\Z,\dotsc,x_n+q_n\Z)\,,
\end{align*}
we obtain the injective homomorphism
\begin{equation}
\Phi\DP B^{\Z_{q_1}\!\times\Z_{q_2}\!\times\dotsm\times\Z_{q_n}}\lto B^{\Z^n},\quad
f\lmto\Phi(f):=f\nach\p\ .
\end{equation}
whose image is the subgroup of \((q_1,q_2,\dotsc,q_n)\)-periodic maps, where a
map $f\DP\Z^n\lto B$ is called \emph{\((q_1,q_2,\dotsc,q_n)\)-periodic}, if
\begin{equation}
f(x+q_ie_i)\,=\,f(x)\ \quad\text{for all $x\in\Z^n$ and all $i\in\{1,2,\dotsc,n\}$.}
\end{equation}
This injection also goes well together with the difference operators $\underset{^i}\Delta$
on $B^{\Z_{q_1}\!\times\Z_{q_2}\!\times\dotsm\times\Z_{q_n}}$ and $B^{\Z^n}$\!\!. We
have
\begin{equation}
\Phi(\underset{^i}\Delta f)\,=\,\underset{^i}\Delta(\Phi(f))\,,
\end{equation}
because
\begin{align*}
\Phi(\underset{^i}\Delta f)(x)
&\ =\ [(\underset{^i}\Delta f)\nach\p](x)\\
&\ =\ [\underset{^i}\Delta f](\p(x))\\
&\ =\ f(\p(x)+\p(e_i))-f(\p(x))\\[1.8mm]
&\ =\ f(\p(x+e_i))-f(\p(x))\taghere\\[1.8mm]
&\ =\ (f\nach\p)(x+e_i)-(f\nach\p)(x)\\[1.8mm]
&\ =\ [\underset{^i}\Delta(f\nach\p)](x)\\
&\ =\ [\underset{^i}\Delta(\Phi(f))](x)\,.
\end{align*}
Based on this observation, we can identify the maps from
$\Z_{q_1}\!\times\Z_{q_2}\!\times\dotsm\times\Z_{q_n}$ to $B$  with the
\((q_1,q_2,\dotsc,q_n)\)-periodic maps from $\Z^n$ to $B$, and simply write $f$ instead
of $\Phi(f)$. More informally speaking, every \((q_1,q_2,\dotsc,q_n)\)-periodic map from
$\Z^n$ to $B$ gives rise to a well-defined map from
$\Z_{q_1}\!\times\Z_{q_2}\!\times\dotsm\times\Z_{q_n}$ to $B$, and vice versa. One just
has to go from congruence classes to representatives, and vice versa. If, say, a
\((q_1,q_2,\dotsc,q_n)\)-periodic map $f$ from $\Z^n$ to $B$ is given, then the
corresponding map from $\Z_{q_1}\!\times\Z_{q_2}\!\times\dotsm\times\Z_{q_n}$ to
$B$, again denoted $f$\!, is well-defined through
\begin{equation}
  f(x_1+q_1\Z,x_2+q_2\Z,\dotsc,x_n+q_n\Z)\,:=\,f(x_1,x_2,\dotsc,x_n)\,.
\end{equation}
In particular, we may view
$B^{\Z_{q_1}\!\times\Z_{q_2}\!\times\dotsm\times\Z_{q_n}}$ as subset of
$B^{\Z^n}$\!\!\!,
and may search for polyfractal representations of maps of the form
$\Z_{q_1}\!\times\Z_{q_2}\!\times\dotsm\times\Z_{q_n}\lto B$ inside the set of all
previously defined polyfracts. The set of \((q_1,q_2,\dotsc,q_n)\)-periodic polyfracts that
describe maps of that form is the subgroup
\begin{equation}
B\tbinom{X_1,X_2,\dotsc,X_n}{\Z_{q_1},\Z_{q_2},\dotsc,\Z_{q_n}}
\,:=\,B\tbinom{X_1,X_2,\dotsc,X_n}{\Z^n}\,\cap\,B^{\Z_{q_1}\!\times
\Z_{q_2}\!\times\dotsm\times\Z_{q_n}}\,.
\end{equation}

The following corollary is a generalization and, at the same time, a subcase of
Theorem\,\ref{sz.eq}\,. It is a generalization as it does not just deal with the case
$q_1=q_2=\dotsb=q_n=0$. It is a subcase and a corollary as $B^{\Z_{q_1}\!\times
\Z_{q_2}\!\times\dotsm\times\Z_{q_n}}$ is a subset of $B^{\Z^n}\!$ and
Theorem\,\ref{sz.eq} can be applied to the functions in that subset:

\begin{Korollar}\label{cor.eq}
The functions $f\DP\Z_{q_1}\!\times\Z_{q_2}\!\times\dotsm\times\Z_{q_n}\lto B$ with
$\fdeg(f)<\infty$ are exactly the functions that are given by
\((q_1,q_2,\dotsc,q_n)\)-periodic \(B\)-polyfracts, i.e.\
$$
\bigl\{f\in B^{\Z_{q_1}\!\times\Z_{q_2}\!\times\dotsm\times\Z_{q_n}}\!\mit\!\fdeg(f)<\infty\bigr\}
\,=\,B\tbinom{X_1,X_2,\dotsc X_n}{\Z_{q_1},\Z_{q_2},\dotsc,\Z_{q_n}}\,.
$$
Moreover, if $f$ is written as polyfract, having total degree $\deg(f)$ and partial degrees
$\deg_i(f)$, then
$$
\fdeg(f)\,=\,\deg(f)\quad\text{and}\quad\pdeg_i(f)\,=\,\deg_i(f)\quad\text{for $i=1,2,\dotsc,n$.}
$$
\end{Korollar}

From this corollary we know that the maps of finite functional degree, and only those,
can be represented by polyfracts. The polyfractal representation of such maps can then
be very helpful. There are several theorems and lemmas in \cite{aimo} that are almost
obvious for polynomials and become about as obvious in the polyfractal representation.
This yields much shorter proofs and additional insights. One example is
\cite[Th.\,5.2]{aimo} (our inequality \eqref{eq.pfp}), another is \cite[Lem.\,6.2]{aimo}
about the functional degree of tensor products. Here, it is actually a hindering limitation
that \cite[Lem.\,6.2]{aimo} is restricted to integral domains $(B,+,\cdot\,)$ as codomain
of the considered functions. With the polynomial representation, however, it is possible
to see beyond integral domains. The additional generality helps then to settle the
question raised in \cite[Problem\,8.3]{aimo}. We will demonstrate that in the last section,
where we explain our solution to Aichinger and Moosbauer's problem.

\section{Classification of Functions of Finite Functional Degree}\label{sec.clas}

From the last corollary we know that the polyfractal maps, i.e.\ those representable by
polyfracts, are the maps of finite functional degree. But, we still have not seen which
maps these are. 
There are two sources. On one hand, for finite commutative groups $A$ and $B$, the
question which maps in $B^A$ are polyfractal was answered in \cite[Th.\,3.7]{sch}. On
the other hand, the question which maps in $B^A$ have finite functional degree was
answered in \cite[Th.\,9.4]{aimo}. As both questions are equivalent, both answers are
the same. To be able to express the solution in an elegant way, we introduce some
notation. We denote with $p_1,\dotsc,p_t$ be the prime divisors of $\abs{A}\abs{B}$,
and with \randd$"A_j"$ and \randd$"B_j"$ the corresponding \(p_j\)"~primary
components (Sylow \(p_j\)-subgroups) of $A$ and $B$, respectively. With $A_j=\{0\}$ if
$p_j\nmid\abs{A}$ and $B_j=\{0\}$ if $p_j\nmid\abs{B}$, we can write $A$ and $B$ as
Cartesian products of $t$ components each. We can write
\begin{equation}
A\,=\,A_1\times A_2\times\dotsb\times A_t\quad\ \text{and}\quad\
B\,=\,B_1\times B_2\times\dotsb\times B_t\,.
\end{equation}
If we further identify each \(t\)-tuple $(f_1,f_2,\dotsc,f_t)$ of functions $f_j\in B_j^{A_j}$
with the function
\begin{align*}
A_1\times A_2\times\dotsb\times A_t&\lto B_1\times B_2\times\dotsb\times B_t\taghere\\
(x_1,x_2,\dotsc,x_t)&\lmto(f_1(x_1),f_2(x_2),\dotsc,f_n(x_t))
\end{align*}
then the set of those \(t\)-tuple of functions becomes a subset of $B^A$\!\!, i.e.\
\begin{equation}
B_1^{A_1}\times B_2^{A_2}\times\dotsb\times B_t^{A_t}\,\sb\,B^A\,.
\end{equation}

With this notation, Theorem\,3.17 in \cite{sch}, combined with the previous corollary,
yields the following classification:

\begin{Satz}
\label{sz.cp} Let $A=A_1\times A_2\times\dotsb\times A_t$ and $B=B_1\times
B_2\times\dotsb\times B_t$ be finite commutative groups, written as direct product of
their (possibly trivial) \(p_j\)"~primary components, as explained above. Then
$$
\bigl\{f\in B^A\!\mit\!\fdeg(f)<\infty\bigr\}
\,=\,B_1^{A_1}\times B_2^{A_2}\times\dotsb\times B_t^{A_t}
\ .
$$
\end{Satz}

\begin{Korollar}\label{cor.cp}
If $A$ and $B$ are finite commutative \(p\)-groups, to the same prime $p$, then all
maps between $A$ and $B$ have finite functional degree.

If there are primes $p_1\neq p_2$ such that $p_1$ divides $\abs{A}$ and $p_2$ divides
$B$, then there are maps of infinite functional degree between $A$ and $B$.
\end{Korollar}

\begin{Beispiel}
If $A=\Z_{60}$ and $B=\Z_{126}\times\Z_{7}$, then $\abs{A}\abs{B}$ has the four prime
divisors $p_1=2$, $p_2=3$, $p_3=5$ and $p_4=7$,
\begin{equation}
A=\Z_{4}\!\times\Z_{3}\!\times\Z_{5}\!\times\{0\}\quad\text{and}\quad
B=\Z_{2}\!\times\Z_{9}\!\times\{0\}\!\times\Z_{7}^2\,.
\end{equation}
According to Theorem\,\ref{sz.cp}, we need to consider the subset
\begin{equation}
\Z_{2}^{\Z_{4}}\!\times\Z_{9}^{\Z_{3}}\!\times\{0\}^{\Z_{5}}\!\times(\Z_{7}^2)^{\{0\}}
\,\sb\,B^A\,,
\end{equation}
of functions of the form
\begin{equation}
(x_1,x_2,x_3,x_4)\lmto\bigl(f(x_1),g(x_2),0,(c,d)\bigr)
\end{equation}
with $f\in\Z_{2}^{\Z_{4}}$, $g\in\Z_{9}^{\Z_{3}}$ and $(c,d)\in\Z_{7}^2$. This can be
simplified, if we omit the zero components that we needed in order to apply
Theorem\,\ref{sz.cp}. With
\begin{equation}
A=\Z_{4}\!\times\Z_{3}\!\times\Z_{5}\quad\text{and}\quad
B=\Z_{2}\!\times\Z_{9}\!\times\Z_{7}\!\times\Z_{7}\,,
\end{equation}
the functions of finite functional degree between $A$ and $B$ have the general form
\begin{equation}
(x_1,x_2,x_3)\lmto\bigl(f(x_1),g(x_2),c,d\bigr)\,.
\end{equation}
If we want to write such a functions as polyfract, we should first write its components as
polyfracts. One might use Theorem\,\ref{sz.tay2} to do so, if the components are given.
In our example, however, $f,g,c,d$ are not given. Hence, $f(x_1)$ can be any polyfract
over $\Z_2$ of degree less than $4$. As $2$ is prime, \cite[Cor.\,3.11]{sch} says that
those are precisely the polyfracts in $\Z_2\tbinom{x_2}{\Z_4}$. It is more complicated to
describe all possible polyfracts $g(x_2)\in\Z_9\tbinom{x_2}{\Z_3}$, as $9$ is not prime.
But, with \cite[Th.\,3.10]{sch} one can easily check if a given polyfract
$g(x_2)\in\Z_9\tbinom{x_2}{\Z}$ actually is \(3\)-periodic and lies in
$\Z_9\tbinom{x_2}{\Z_3}$. For instance, we could have
\begin{equation}
f(x_1)=\tbinom{x_1}{3}+\tbinom{x_1}{1}\,,\ \ g(x_2)=6\tbinom{x_2}{1}+3\,,\ \ c=4\,,\ \ d=5\,.
\end{equation}
These choices can now be combined to the following example of a polyfractal map from
$A$ to $B$:
\begin{align*}\label{eq.ex}
\Z_{4}\!\times\Z_{3}\!\times\Z_{5}\,\lto\,
&\,\Z_{2}\!\times\Z_{9}\!\times\Z_{7}\!\times\Z_{7}\\
(x_1,x_2,x_3)\,\lmto\,
&\biggl(\dbinom{x_1}{3}+\dbinom{x_1}{1}\,,\,6\dbinom{x_2}{1}+3\,,\,4\,,\,5\biggr)\\[5pt]
&=\,(1,0,0,0)\dbinom{x_1,x_2,x_3}{3\,,\,0\,,\,0}+(1,0,0,0)\dbinom{x_1,x_2,x_3}{1\,,\,0\,,\,0}\taghere\\
&\quad\,+(0,6,0,0)\dbinom{x_1,x_2,x_3}{0\,,\,1\,,\,0}+(0,3,4,5)\\[-5mm]
\end{align*}
\end{Beispiel}

In the last representation in \eqref{eq.ex}, one may observe how Theorem\,\ref{sz.cp}
ensures that the coefficients in $\Z_{2}\!\times\Z_{9}\!\times\Z_{7}\!\times\Z_{7}$ have
many zeros. In general, we have the following:

\begin{Korollar}
Let $p_1,p_2,\dotsc,p_t$ be different primes. For $j=1,2,\dotsc,t$, let $B_{j}$ be a finite
commutative \(p_j\)-group, and let $q_{j,1},\,q_{j,2},\,\dotsc,\,q_{j,n_j}$ be powers of
$p_j$ (where $B_{j}$ and some $q_{j,i}$ may be trivial but $n_j\geq1$).
If a non-constant term
$$(b_1,\dotsc,b_t)
\dbinom{X_{1,1},\dotsc,X_{1,n_1},\,X_{2,1},\dotsc,X_{2,n_2},\,\dotsc,\,X_{t,1},\dotsc,X_{t,n_t}}%
{\d_{1,1},\dotsc,\,\d_{1,n_1}\,,\ \d_{2,1},\dotsc,\,\d_{2,n_2}\,,\ \dotsc\,,\d_{t,1},\dotsc,\d_{t,n_t}}$$ %
occurs in a
\((q_{1,1},\dotsc,q_{1,n_1},q_{2,1},\dotsc,q_{2,n_2},\dotsc,q_{t,1},\dotsc,q_{t,n_t})\)-periodic
polyfract, i.e.\ in the standard expansion of a polyfract in
$$(B_1\times B_2\times\dotsb\times B_t)
\dbinom{X_{1,1},\dotsc,X_{1,n_1}\,,\,X_{2,1},\dotsc,X_{2,n_2}\,,\,\dotsc,\,\,X_{t,1},\dotsc,X_{t,n_t}\,}%
{\Z_{q_{1,1}},\dotsc,\Z_{q_{1,n_1}},\,
\Z_{q_{2,1}},\dotsc,\Z_{q_{2,n_2}},\,\dotsc,\,\Z_{q_{t,1}},\dotsc,\Z_{q_{t,n_t}}}\,,$$ then
there exists a $j_0\in\{1,2,\dotsc,t\}$ such that
$$\d_{j,1},\,\d_{j,2},\,\dotsc,\,\d_{j,n_j}=0\ \ \text{and}\ \ b_j=0\ \ \text{for all $j\neq j_0$.}$$
\end{Korollar}

The corollary applies to all polyfractal maps between arbitrary nontrivial finite
commutative groups $A$ and $B$. One can choose $p_1,p_2,\dotsc,p_t$ as the prime
divisors of $\abs{A}\abs{B}$, and $A_j:=\Z_{q_{j,1}},\dotsc,\Z_{q_{j,n_j}}$ and $B_j$ as
\(p_j\)-primary component of $A$ and $B$, respectively. If $A_j$ is trivial, one just has
to select $n_j$ as $1$ and $(q_{j,1},\,q_{j,2},\,\dotsc,\,q_{j,n_j})$ as $(p^0)$, and treat
$A_j=\{0\}$ as $\Z_1$. If a map $f\DP A\lto B$ has finite functional degree, then $f\in
B_1^{A_1}\times B_2^{A_2}\times\dotsb\times B_t^{A_t}$ and we can write it as
$f=(f_1,f_2,\dotsc,f_t)$ with components $f_j\in B_j^{A_j}$. Those components $f_j$ are
easily read off if $f$ is written as polyfract, as in the corollary. Each non-constant term
\begin{equation}
(b_1,\dotsc,b_t)
\dbinom{X_{1,1},\dotsc,X_{1,n_1},\,X_{2,1},\dotsc,X_{2,n_2},\,\dotsc,\,X_{t,1},\dotsc,X_{t,n_t}}%
{\d_{1,1},\dotsc,\,\d_{1,n_1}\,,\ \d_{2,1},\dotsc,\,\d_{2,n_2}\,,\ \dotsc\,,\ d_{t,1},\dotsc,\d_{t,n_t}}
\end{equation}
of $f$ contributes the term
\begin{equation}
b_{j_0}
\dbinom{X_{j_0,1},\dotsc,X_{j_0,n_{j_0}}}{\d_{j_0,1},\dotsc,\,\d_{j_0,n_{j_0}}}
\end{equation}
to the polyfractal expansion of $f_{j_0}$, where $j_0$ is the unique index with
$b_{j_0}\neq0$. The constant term is the only term of $f$ that may contribute to more
than one $f_j$. Each $f_j$ simply gets the \(j^{\text{th}}\) entry of that term. It is as in
\eqref{eq.ex}, in the example above. In particular, we see that
\begin{equation}\label{eq.max2}
\deg(f)=\max\bigl(\deg(f_1),\deg(f_2),\dotsc,\deg(f_t)\bigr)\,,
\end{equation}
where the degree of a function $f_j$ with trivial domain is $0$. This is because the trivial
domain is represented as $\Z_1$, so that $f_j$, seen as function on $\Z$, is
\(1\)-periodic, i.e.\ constant.


\section{The Largest Finite Degree}\label{sec.max}

In this section, in Theorem\,\ref{sz.mda} below, we determine the maximal finite degree
of a function between finite commutative groups, i.e.\ the largest degree that a function
of finite degree can have. In particular, this solves (in Equation\,\eqref{eq.sol}) a
problem raised by Aichinger and Moosbauer in \cite[Problem\,8.3]{aimo}, the problem of
finding the nilpotentcy degree of the augmentation ideal of the group ring
$\Z_{p^{\b_1}}[\Z_{p^{\ä_1}}\!\times\Z_{p^{\ä_2}} \!\times\dotsm\times\Z_{p^{\ä_n}}]$,
where $p$ is prime. Actually, if one is aware of the connection between polyfracts and
functions of finite functional degree, our Theorem\,3.9 in \cite{sch} already solved this
more special problem.
We repeat the steps that lead to this solution and involve polyfracts, to explain the role
of the polyfractal representation, before we further generalize things. In doing so, we
also point to some simplifications in proofs of underlying facts, based on a paper that we
were not aware of when we wrote \cite{sch}. We start the investigations with very
special types of commutative groups, and then generalize the results from subsection to
subsection.

\subsection{Between \,$\Z_{p^\ä}$ and \,$\Z_{p^\b}$}

Let $\ä,\b,p\in\Z^+$ with $p$ prime. The Lagrange function (characteristic function)
$\chi\DP\Z_{p^\ä}\!\lto\Z_{p^\b}$ is defined by
\begin{equation}
\chi(x):=\begin{cases}
1
&\text{if $x=0$,}\\
0
&\text{otherwise.}
\end{cases}
\end{equation}
In this subsection, we want to explain why
\begin{equation}\label{eq.nt1}
\fdeg(\chi)\,=\,p^{\ä}-1+(\b{-}1)(p{-}1)p^{\ä-1}\,=\,\b p^\ä-(\b{-}1)p^{\ä-1}-1\,.
\end{equation}

Actually, we have deduced this in \cite[Th.\,3.9]{sch} already, but the proof there is not
short and relies to one half on \cite[Cor.\,4.16]{hr}. As we are considering only the cyclic
case in Equation\,\eqref{eq.nt1}, however, the result \cite[Cor.\,4.16]{hr} alone is strong
enough to deduce Equation\,\eqref{eq.nt1}, without the need for polyfractal
representations. If we look into \cite{hr}, we see that the largest possible functional
degree can be found by determining the isomorphy type of the \(\Z\)-module
$\Phi_k(\Z_{p^\ä},\Z_{p^\b})$ of all functions $f\DP\Z_{p^\ä}\!\to\Z_{p^\b}$ of functional
degree at most $k$. If that isomorphy type is known, it can simple be observed that
\begin{equation}
\Abs{\Phi_{m-1}(\Z_{p^\ä},\Z_{p^\b})}<\Abs{\Phi_m(\Z_{p^\ä},\Z_{p^\b})}=\Abs{\Z_{p^\b}^{\Z_{p^\ä}}}
\,\ \text{if}\,\ m=\b p^\ä-(\b{-}1)p^{\ä-1}-1\,.
\end{equation}
This shows that there exists a function $f\DP\Z_{p^\ä}\!\to\Z_{p^\b}$ of degree $m$, but
none of degree bigger than $m$. It follows that $\chi$ must have that degree $m=\b
p^\ä-(\b{-}1)p^{\ä-1}-1$. That is what we want, but the proof is long and quite indirect.

Recently, we discovered in \cite[Th.\,10]{wi} a more direct proof. In this proof one
considers the Taylor expansion of $\chi$, as in our Theorem\,\ref{sz.tay2} with $n=1$.
It is easy to calculate the coefficient $c_{\d,\ä}$ of $\tbinom{x}{\d}$ in that expansion
(see e.g.\ \cite[Th.\,3.6]{sch} or \cite[proof of Th.\,10]{wi}). A representative $\hat
c_{\d,\ä}\in\Z$ of the congruence class $c_{\d,\ä}\in\Z_{p^\b}$ is given by
\begin{equation}
\hat c_{\d,\ä}\,=\,\sum_i(-1)^i\dbinom{\d}{i}\,\,\ ,
\end{equation}
where the sum is taken over all $i\in\{0,1,\dotsc,\d\}$ with $i\equiv\d\pmod{p^\ä}$. The
tricky part is then to prove the following lemma:

\begin{Lemma}\label{lem.dpb}
Let $\ä,\b,\d\in\Z^+$.\smallskip

If $\d>\b p^\ä-(\b{-}1)p^{\ä-1}-1$, then $p^\b\div\hat c_{\d,\ä}$.\smallskip

If $\d=\b p^\ä-(\b{-}1)p^{\ä-1}-1$, then $p^\b\ndiv\hat c_{\d,\ä}$.
\end{Lemma}

\begin{Beweis}
We use the fact that $\hat c_{\d,\ä}$ is also equal to the coefficient of $x^{\d\bmod
p^\ä}$ when $(x-1)^\d$ is reduced modulo $x^{p^\ä}\!\!-1$. For
\begin{equation}
\d=\b
p^\ä-(\b{-}1)p^{\ä-1}-1=(\b(p-1)+1)p^{\ä-1}-1
\end{equation}
this reduction is provided in \cite[Eq.\,(22)]{wi}, at least modulo $p^\b$\!. It says that
\begin{equation}\label{eq.red}
(x-1)^\d\equiv(-p)^{\b-1}\sum_{j=0}^{p^\ä-1}x^j\pmod{p^\b,x^{p^\ä}\!\!-1}\,.
\end{equation}
So, $\hat c_{\d,\ä}\equiv(-p)^{\b-1}\not\equiv0\pmod{p^\b}$ if $\d=\b
p^\ä-(\b{-}1)p^{\ä-1}-1$. If we multiply \eqref{eq.red} with $x-1$ then the right side
becomes zero, if reduced modulo $x^{p^\ä}\!\!-1$. Hence, $\hat
c_{\d,\ä}\equiv0\pmod{p^\b}$ whenever $\d>\b p^\ä-(\b{-}1)p^{\ä-1}-1$.
\end{Beweis}

Based on \cite[Cor.\,4.16]{hr}, one may also prove this lemma as in \cite[Cor.\,3.5\ \&\
Eq.\,(88)]{sch}. It tells us that, in the Taylor expansion of
$\chi\DP\Z_{p^\ä}\!\to\Z_{p^\b}$, the coefficient $c_{\d,\ä}$ of $\tbinom{x}{\d}$ is
nonzero for $\d=\b p^\ä-(\b{-}1)p^{\ä-1}-1$ but zero if $\d>\b p^\ä-(\b{-}1)p^{\ä-1}-1$. So,
the degree of $\chi$ is $\b p^\ä-(\b{-}1)p^{\ä-1}-1$, which shows that
Equation\,\eqref{eq.nt1} holds, indeed.

\subsection{Between $\Z_{p^{\ä_1}}\!\times\Z_{p^{\ä_2}}
\!\times\dotsm\times\Z_{p^{\ä_n}}$ and \,$\Z_{p^\b}$}

The power of Lemma\,\ref{lem.dpb} is not limited to the case of cyclic \(p\)-groups. If we
work with polyfractal representations, it can easily tell us what happens on the more
general domain $\Z_{p^{\ä_1}}\!\times\Z_{p^{\ä_2}} \!\times\dotsm\times\Z_{p^{\ä_n}}$.
We define the multivariate Lagrange function
 $\chi\DP\Z_{p^{\ä_1}}\!\times\Z_{p^{\ä_2}}\!\times\dotsm\times\Z_{p^{\ä_n}}\lto\Z_{p^\b}$
as product of the $n$ univariate Lagrange functions
$\chi_j\DP\Z_{p^{\ä_j}}\lto\Z_{p^\b}$, i.e.
\begin{equation}
\chi(x_1,x_2,\dotsc,x_n)\,=\,\chi_1(x_1)\,\chi_2(x_2)\dotsm\chi_n(x_n)\,.
\end{equation}
We want to show that
\begin{equation}\label{eq.chi}
\fdeg(\chi)\,=\,\sum_{j=1}^np^{\ä_j}-n+(\b{-}1)(p{-}1)p^{\ä_{\text{max}}-1}\,,
\end{equation}
where $\ä_{\text{max}}:=\max\limits_{1\leq j\leq n}\ä_j$.

The coefficient of $\tbinom{x_j}{\d_j}$ in the expansion of $\chi_j(x_j)$ is $c_{\d_j,\ä_j}$.
Hence, the coefficient of the monofract
$\tbinom{x_1,x_2,\dotsc,x_n}{\d_1,\d_2,\dotsc,\d_n}$ in the product
$\chi_1(x_1)\,\chi_2(x_2)\dotsm\chi_n(x_n)$ is $c_{\d_1,\ä_1}c_{\d_2,\ä_2}\dotsm
c_{\d_n,\ä_n}=\hat c_{\d_1,\ä_1}\hat c_{\d_2,\ä_2}\dotsm\hat c_{\d_n,\ä_n}+p^\b\Z$,
and with Lemma\,\ref{lem.dpb} it is easy to see when that coefficient is zero. We may
assume $\ä_1=\ä_{\text{max}}$. If
\begin{equation}\label{eq.max}
\d_1=\b p^{\ä_1}-(\b{-}1)p^{\ä_1-1}-1\ ,\,\ \d_2=p^{\ä_2}-1\ ,\ \dotsc\ ,\,\ \d_n=p^{\ä_n}-1\ ,
\end{equation}
then Lemma\,\ref{lem.dpb} tells us that
\begin{equation}
p^{\b}\ndiv\hat c_{\d_1,\ä_1}\ ,\quad
p\ndiv\hat c_{\d_2,\ä_2}\ ,\ \dotsc\ ,\quad
p\ndiv\hat c_{\d_n,\ä_n}\,\ ,
\end{equation}%
so that $c_{\d_1,\ä_1}c_{\d_2,\ä_2}\dotsm c_{\d_n,\ä_n}\neq0\in\Z_{p^\b}$. This shows
that
\begin{equation}\label{eq.chi2}
\fdeg(\chi)\,\geq\,\sum_{j=1}^np^{\ä_j}-n+(\b{-}1)(p{-}1)p^{\ä_{\text{max}}-1}\,.
\end{equation}

To examine more precisely from which total degree $\d_1+\d_2+\dotsb+\d_n$ onwards
the coefficients $c_{\d_1,\ä_1}c_{\d_2,\ä_2}\dotsm c_{\d_n,\ä_n}$ start to vanish in
$\Z_{p^\b}$, we observe what happens if we increase the different $\d_j$. Roughly
speaking, Lemma\,\ref{lem.dpb} tells us that the bigger we make the $\d_j$, the more
often $p$ will divide $\hat c_{\d_1,\ä_1}\hat c_{\d_2,\ä_2}\dotsm\hat c_{\d_n,\ä_n}$.
Looking at one fixed index $j$, and starting with $\d_j=0$, we see that the first
$p^{\ä_j}-1$ increments of $\d_j$ are for free. As long as $\d_j\leq p^{\ä_j}-1$, we have
$p\ndiv\hat c_{\d_j,\ä_j}$. After reaching the first threshold $p^{\ä_j}-1$, the next
increment of $\d_j$ will make $p$ divide $\hat c_{\d_j,\ä_j}$, but we can increase $\d_j$
by up to $p^{\ä_j}-p^{\ä_j-1}$ before we reach the next threshold.
Each further increase of $\d_j$ 
by $p^{\ä_j}-p^{\ä_j-1}$ can be bought for just one additional $p$ as divisor of $\hat
c_{\d_j,\ä_j}$. Moreover, regarding increments beyond the first threshold, increasing
$\d_1$ is
at least as cheap as increasing any other $\d_j$, because 
$p^{\ä_1}-p^{\ä_1-1}\geq p^{\ä_j}-p^{\ä_j-1}$.
This is why we obtain the largest possible sum $\d_1+\d_2+\dotsb+\d_n$ for which still
$p^\b\ndiv\hat c_{\d_1,\ä_1}\hat c_{\d_2,\ä_2}\dotsm\hat c_{\d_n,\ä_n}$ by just
increasing $\d_1$ beyond its first threshold $p^{\ä_1}-1$ to
$p^{\ä_1}-1+(\b-1)(p^{\ä_1}-p^{\ä_1-1})=\b p^{\ä_1}-(\b{-}1)p^{\ä_1-1}-1$, and leaving all
the other $\d_j$ at their first threshold, as in \eqref{eq.max}. Any coefficient
$c_{\d_1,\ä_1}c_{\d_2,\ä_2}\dotsm c_{\d_n,\ä_n}$ with larger sum
$\d_1+\d_2+\dotsb+\d_n$ is zero in $\Z_{p^\b}$. So, \eqref{eq.chi2} is best possible
and \eqref{eq.chi} holds, indeed.
%

Since every function
$f\DP\Z_{p^{\ä_1}}\!\times\Z_{p^{\ä_2}}\!\times\dotsm\times\Z_{p^{\ä_n}}\lto\Z_{p^\b}$ is
a linear combination of shifted Lagrange functions $\chi(x-a)$,
the degree of $f$ is bounded by the degree of $\chi$:


\begin{Satz}\label{sz.md}
Let $p$ be a prime and $\ä_1,\ä_2,\dotsc,\ä_n,\b\geq1$. Maps of the form
\begin{equation*}
f\DP\Z_{p^{\ä_1}}\!\times\Z_{p^{\ä_2}}
\!\times\dotsm\times\Z_{p^{\ä_n}}\lto\, \,\Z_{p^{\b}}
\end{equation*}
have bounded functional degree, and the best upper bound on the functional degree is
given by
$$
\fdeg(f)\,\leq\,\sum_{j=1}^np^{\ä_j}-n+(\b{-}1)(p{-}1)p^{\ä_{\text{max}}-1}\,,
$$
where $\ä_{\text{max}}:=\max\limits_{1\leq j\leq n}\ä_j$.
\end{Satz}

Aichinger and Moosbauer explained in \cite[Lemma\,7.3]{aimo} that, for functions
$f\DP\Z_{p^{\ä_1}}\!\times\Z_{p^{\ä_2}} \!\times\dotsm\times\Z_{p^{\ä_n}}\lto\,
\,\Z_{p^{\b}}$, the best upper bound of the degree plus one is the nilpotentcy degree
$\nu$ of the augmentation ideal of the group ring
$\Z_{p^\b}[\Z_{p^{\ä_1}}\!\times\Z_{p^{\ä_2}} \!\times\dotsm\times\Z_{p^{\ä_n}}]$.
Finding this degree was posted as open problem in \cite[Problem\,8.3]{aimo}. The
solution is
\begin{equation}\label{eq.sol}
\nu\,=\,\sum_{j=1}^np^{\ä_j}-n+1+(\b{-}1)(p{-}1)p^{\ä_{\text{max}}-1}\,.
\end{equation}

\subsection{Between \,$\Z_{p^{\ä_1}}\!\times\Z_{p^{\ä_2}}
\!\times\dotsm\times\Z_{p^{\ä_n}}$ and \,$\Z_{p^{\b_1}}\!\times\Z_{p^{\b_2}}
\!\times\dotsm\times\Z_{p^{\b_t}}$}

Maps of the form
\begin{equation}
f\DP\Z_{p^{\ä_1}}\!\times\Z_{p^{\ä_2}}
\!\times\dotsm\times\Z_{p^{\ä_n}}\lto\, \,\Z_{p^{\b_1}}\!\times\Z_{p^{\b_2}}
\!\times\dotsm\times\Z_{p^{\b_t}}
\end{equation}
can be split into $t$ maps
$f_1,f_2,\dots,f_t$ such that
\begin{equation}
f(x)\,=\,\bigl(f_1(x),f_2(x),\dotsc,f_t(x)\bigr)\,.
\end{equation}
In that case, simply
\begin{equation}\label{eq.t1}
\fdeg(f)\,=\,\max_{1\leq i\leq t}(\fdeg(f_i))\,.
\end{equation}
This is obvious in polyfractal representation, but also easy to see from the definition of
the functional degree (see \cite[Lem.\,3.4]{aimo}). Hence,
the best upper bound on the functional degree is given by
\begin{equation}\label{eq.pgr}
\fdeg(f)\,\leq\,\sum_{j=1}^np^{\ä_j}-n+(\b_{\text{max}}{-}1)(p{-}1)p^{\ä_{\text{max}}-1}\,,
\end{equation}
where $\ä_{\text{max}}:=\max\limits_{1\leq j\leq n}\ä_j$ and
$\b_{\text{max}}:=\max\limits_{1\leq i\leq t}\b_i$.\smallskip

\subsection{Between Arbitrary Finite Commutative Groups}

If we increase the generality further, to maps $f\DP A\lto B$ between arbitrary finite
commutative groups, then the functional degree may become infinite. Assume,
$p_1,\dotsc,p_t$ are the prime divisors of $\abs{A}\abs{B}$, with corresponding
\(p_i\)"~primary components $A_i$ and $B_i$ of $A$ and $B$, respectively. If $t=1$,
everything is finite, by Corollary\,\ref{cor.cp}\,. But, if $t>1$ and $\abs{A}>1$ and
$\abs{B}>1$, then there exist a $j_1$ and a $j_2\neq j_1$ with $\abs{A_{j_1}}>1$ and
$\abs{B_{j_2}}>1$. In this situation, functions of infinite degree are easy to find. Every
map that sends a non-trivial elements of $A_{j_1}$ to a non-trivial element of $B_{j_2}$
does not lie in $\prod_{j=1}^t{B_j^{A_j}}$ and must have infinite degree, by
Theorem\,\ref{sz.cp}.

If we restrict ourselves to functions of finite degree, then $f$ splits again, but differently.
It does not split as in the last subsection. In this case, by Theorem\,\ref{sz.cp},
$f\in\prod_{j=1}^t{B_j^{A_j}}$ and there are maps $f_i\DP A_i\lto B_i$ such that
\begin{equation}
f(x_1,x_2,\dotsc,x_n)\,=\,(f_1(x_1),f_2(x_2),\dotsc,f_t(x_n))\,.
\end{equation}
Thus, by \eqref{eq.max2},
\begin{equation}
\fdeg(f)=\max\bigl(\fdeg(f_1),\fdeg(f_2),\dotsc,\fdeg(f_t)\bigr)\,.
\end{equation}
Since the best upper bound for $\fdeg(f_j)$ is given in \eqref{eq.pgr} if $A_j\neq\{0\}$
and $B_j\neq\{0\}$, and $\fdeg(f_j)=0$ if $A_j=\{0\}$ or $B_j=\{0\}$, we can easily
calculate the best upper bound for $\fdeg(f)$. Of course, if $\gcd(\abs{A},\abs{B})=1$
then the only functions in $\prod_{j=1}^t{B_j^{A_j}}$
are the constant functions, which have degree $0$. So, things only get interesting if
$\gcd(\abs{A},\abs{B})>1$. We obtain the following result:

\begin{Satz}\label{sz.mda}
Let $A$ and $B$ be finite commutative groups of non-coprime order, and let
$p_1,p_2,\dotsc,p_s$ be the prime divisors of $\gcd(\abs{A},\abs{B})>1$. For
$i=1,2,\dotsc,s$, denote with $A_i$ resp.\,$B_i$ the (non-trivial) \(p_i\)"~primary
component of $A$ resp.\,$B$. Assume
\begin{align*}
A_i&=\Z_{q_{i,1}}\!\times\dotsm\times\Z_{q_{i,n_i}}\!\!&\text{with}\,\ q_{i,j}&=p_i^{\ä_{i,j}}>1
\,,&\!\text{and set}\,\ \ä_{i,\text{max}}&:=\max_{1\leq j\leq n_i}\!(\ä_{i,j})\,,\\
B_i&=\Z_{r_{i,1}}\!\times\dotsm\times\Z_{r_{i,m_i}}\!\!&\text{with}\,\ r_{i,j}&=p_i^{\b_{i,j}}>1
\,,&\!\text{and set}\,\ \b_{i,\text{max}}&:=\max_{1\leq j\leq m_i}\!(\b_{i,j})\,.
\end{align*}
For functions $f\DP A\lto B$ of finite functional degree, the best upper bound on the
functional degree is given by
$$
\fdeg(f)\,\leq\,\max_{1\leq i\leq s}\Bigl(\,\sum_{j=1}^{n_i}p_i^{\ä_{i,j}}-n_i+
(\b_{i,\text{max}}{-}1)(p_i{-}1)p_i^{\ä_{i,\text{max}}-1}\Bigr)\,.\medskip
$$
\end{Satz}



\begin{thebibliography}{00}
\bibitem
{aimo}
  E.\, Aichinger, J.\,Moosbauer: Chevalley-Warning type results on abelian groups.
  \textit{Journal of Algebra, Volume 569 (2021), 30-66.}
\bibitem
{hr}
  J.\,Hrycaj: Polynomial Mappings of Finitely Generated \(\Z\)"~Modules.\\
  \textit{Journal of Number Theory 26(3) (1987), 308-324.}
\bibitem
{la}
  M.\,Laczkovich:
  Polynomial Mappings on Abelian Groups.\\
  \textit{Aequationes Math.\ 68 (2004), 177-199.}
\bibitem
{ma}
   P.\,Mayr: The Subpower Membership Problem for Mal'cev Algebras.\\
  \textit{International Journal of Algebra and Computation, 22(7) (2012), doi:10.1142/S0218196712500750.}
\bibitem
{sch}
  U.\,Schauz: Classification of polynomial mappings between commutative groups.
  \textit{Journal of Number Theory, Volume 139 (2014), 1-28.}
\bibitem
{wi} M.\,Wilson: A Lemma on Polynomials Modulo $p^m$ and Applications to Coding
Theory. \textit{Discrete Mathematics 306 (23) (2006), 3154-3165}
\end{thebibliography}
\end{document}